\documentclass[11pt,twoside]{article}

\usepackage{epsfig,amsfonts,xcolor}
\usepackage{amsmath}
\usepackage{lineno}
\usepackage{pdflscape}
\usepackage{caption}  
\usepackage{subcaption}
\usepackage{threeparttable} 
\usepackage{multirow}

\usepackage{amssymb, palatino,geometry,url}
\usepackage{algorithmic}
\usepackage[colorlinks=true,linkcolor=blue,citecolor=blue,urlcolor=blue,breaklinks=true]{hyperref}

\usepackage[authoryear,round]{natbib}
\setcitestyle{authoryear,round,aysep={,},yysep={,},citesep={;},maxcitenames=1}

\usepackage{tikz}
\usetikzlibrary{shapes.geometric}
\DeclareRobustCommand{\legendnode}[2]{%
  \tikz[baseline=-0.6ex]\node[draw, fill=#1, #2, minimum size=6pt, inner sep=0pt] {};
}
\usepackage[dvipsnames]{xcolor}
\definecolor{mskyblue}{HTML}{87CEEB}       
\definecolor{mseagreen}{HTML}{3CB371}      
\definecolor{mgoldenrod}{HTML}{DAA520}     
\definecolor{mlightcoral}{HTML}{F08080}    
\DeclareUnicodeCharacter{2212}{\ensuremath{-}}

\allowdisplaybreaks

\usepackage{fancyhdr}
\newcommand{\be}{\begin{equation}}
\newcommand{\ee}{\end{equation}}
\newcommand{\bea}{\begin{eqnarray}}
\newcommand{\eea}{\end{eqnarray}}

\newcommand{\bvec}{\left(\begin{array}{c}}
\newcommand{\evec}{\end{array}\right)}
\newcommand{\bsub}{\begin{subequations}}
\newcommand{\esub}{\end{subequations}}

\begin{document}

\title{\Large Simulation-Based Optimization over Discrete Spaces using Projection to Continuous Latent Spaces}

\author{Gabriel Hernández-Morales${}^{\dag,\ddag}$,  Brenda Cansino-Loeza${}^{\dag}$\\ Arturo Jiménez-Gutiérrez${}^{\ddag}$,  Victor M. Zavala${}^{\dag}$\thanks{Corresponding Author: victor.zavala@wisc.edu.}\\
 \\
  {\small ${}^{\dag}$Department of Chemical and Biological Engineering}\\
 {\small \;University of Wisconsin - Madison, 1415 Engineering Dr, Madison, WI 53706, USA}\\
 {\small ${}^{\ddag}$ Department of Chemical Engineering}\\
 {\small \; Tecnologico Nacional de Mexico - Instituto Tecnologico de Celaya,  Celaya,  Guanajuato, Mexico}}
 \date{}
\maketitle

\begin{abstract}

Simulation-based optimization of complex systems over discrete decision spaces is a challenging computational problem. Specifically, discrete decision spaces lead to a combinatorial explosion of possible alternatives, making it computationally prohibitive to perform simulations for all possible combinations. In this work, we present a new approach to handle these issues by transforming/projecting the discrete decision space into a continuous latent space using a probabilistic model know as Variational AutoEncoders. The transformation of the decision space facilitates the implementation of Bayesian optimization (BO), which is an efficient approach that strategically navigates the space to reduce the number of expensive simulations. Here, the key observation is that points in the latent space correspond to decisions in the original mixed-discrete space, but the latent space is much easier to navigate using BO. We illustrate the benefits of our approach through a couple of case studies that aim to design complex distillation systems: the recovery of caprylic acid from water by liquid-liquid extraction and the separation of an azeotropic mixture using a thermally couple column know as extractive dividing wall column. 

\end{abstract}


\section{Introduction}

The design of complex systems, such as chemical processes, often involves a search space over discrete variables (e.g., topology of the chemical process, solvent selection, or equipment selection). The presence of discrete decisions complicates the use of optimization techniques, particularly when the system model is a black-box simulator (there is no explicit access to algebraic equations). This challenge has motivated the development of black-box optimization strategies, which rely on observed input–output data obtained from the simulator  \citep{biegler2018new, rogers2015feasibility}. Among these approaches, genetic algorithms (GAs) and other evolutionary algorithms have been widely used and applied to various processes with different levels of complexity \citep{ZHU2025126703, PANDIT2022121437, JING2025131361, Ibrahim2017}. Evolutionary algorithms are meta-heuristic approaches that can handle discrete design spaces, but they might require a large number of simulations and there are no clear guidelines on when to stop the search. Similar limitations are often encountered with other meta-heuristic approaches \citep{JAVALOYESANTON2022107655}. 

Bayesian Optimization (BO) has emerged as a powerful black-box optimization method due to its ability to efficiently explore design spaces using probabilistic surrogate models (usually a Gaussian Process-GP) that are  trained on limited data, making it useful when simulations or experiments are costly. A key aspect of BO is that it leverages uncertainty information of the surrogate model to help explore the design space \citep{paulson2024bayesian, scyphers2024bayesian}. The application of BO to the design of chemical processes, such as the design of separation systems, has shown that this approach is effective at reducing the number of simulations. However, existing studies have focused on decision spaces that only involve continuous decisions, such as distillate-to-feed and reflux ratios \citep{PEREZONES2024109708, PANALESPEREZ2025110349}. In order to handle discrete variables in BO, some studies have used rounding strategies \citep{garrido2020dealing}; this approach, however, can lead the algorithm to get trapped in local optima or can lead to infeasible solutions. For example, the work reported in \citet{Phuc1007} proposed an approach tailored for discrete search spaces that avoids redundant sampling and local optima by optimizing the exploration-exploitation factor and the scaling length of the covariance function. While this method improves performance in purely discrete spaces, the proposed modifications are ad-hoc. The work in \citet{wan2021thinkglobalactlocal}  proposed a BO method for categorical and mixed-discrete domains, combining customized kernels with a local trust region strategy for high-dimensional optimization. However, it is in general difficult to fit a smooth surrogate model over a discrete space. 

Researchers in the molecular design community have proposed to transform the molecular design space (which is inherently discrete) into a continuous space  by using \textit{Variational AutoEncoders} (VAE) \citep{Benjamin2018}. VAEs aim to identify latent representations of input data; this is done by training a neural network (known as an encoder) that compresses the input into a lower-dimensional space and by simultaneously training another neural network (known as a decoder) that aims to reconstruct the original input data from the latent representation. This approach has been widely used to navigate the design space of molecules and proteins. 
Building on this concept, the work in \citet{stanton2022} introduced LaMBO, a BO strategy that employs a low-dimensional latent space and an autoregressive model to optimize small molecules and protein sequences. This methodology was extended with LaMBO-2 in \citet{gruver2023}, which replaces Gaussian process (GP) models with ensemble methods, maintaining the focus on discrete spaces. Other works have also explored the use of BO over discrete inputs, using string-based kernels or exploiting the structure of the latent space in molecular applications \citep{moss2020, Bombarelli7b00572, tripp2020}. More recently, \citet{michael2024} proposed CoRel, an approach for protein sequence optimization that reformulates the discrete BO problem into a continuous one via a VAE latent space. This methodology enables the direct definition of a GP model over probability distributions, utilizing a new kernel based on the weighted Hellinger distance. This formulation not only incorporates prior domain knowledge but also proves to be effective in low-budget optimization scenarios, such as protein design, where data are scarce and experimental evaluations are expensive. The use of autoencoders to transform complex spaces into tractable ones (e.g., from nonlinear to linear) has also found many uses in the development of dynamical models; specifically, diverse studies have reported on the ability of using autoencoders to help discover linear latent spaces under which it is much easier to develop dynamical models and control strategies \citep{schulze2022identification}. 

In this work,  we aim to investigate whether discrete spaces encountered in simulation-based process design can be represented as continuous spaces using VAEs. We use a BO algorithm to navigate the continuous latent space, to guide the strategic collection of simulation data over the design space, and to identify optimal designs. We apply our framework to a couple of case studies arising in the design of complex separation systems. Our results provide interesting insights into the ability of using VAEs to transform complex decision spaces and to facilitate their navigation. Specifically, we show that the proposed approach can help identify designs of high quality with a few simulations. While our results are computational and empirical in nature, we believe that the use of VAEs offer an interesting and practical alternative to help navigate complex decision spaces arising in optimization.  All data and code needed to reproduce our results can be found at \url{https://github.com/zavalab/ML/tree/master/bayesianvae}.

\section{Projection of Discrete into Continuous Spaces}

A variational autoencoder (VAE) is a machine learning model that aims to learn a latent space from which original data can be reconstructed with minimal information loss \citep{Benjamin2018}. A VAE is a compression-decompression system that operates as follows: the {\em encoder} takes points that live in an input discrete space $\mathbf{x} \in \mathcal{X}$ and compresses them using a sequence of nonlinear transformations (conducted using a neural network) into points that live in a continuous latent space $\mathbf{z} \in \mathcal{Z}$. The {\em decoder} uses a sequence of nonlinear transformations (conducted using another neural network) to reconstruct the original points from points in the latent space in a way that it minimizes a measure of reconstruction error. The parameters of the VAE model are trained using available data (in an unsupervised manner). For simplicity in the discussion, we assume that the input space $\mathcal{X}$ is completely discrete (it does not include continuous variables); any continuous variables can be processed via discretization. In practice (e.g., in experiments), one often implements continuous variables in discrete values. The concepts discussed can be generalized to account for input spaces that are mixed-discrete.

Unlike a standard autoencoder, which builds a latent space that is unstructured \citep{Bombarelli7b00572}, the VAE architecture is unique in that it induces a  probabilistic structure into the latent space. Specifically, each data point in the original space is encoded not as a single point but as a probability distribution (typically a Gaussian distribution). During training, samples from this distribution are drawn to force the decoder to reconstruct the original input from noisy/perturbed latent vectors. This probabilistic step ensures that the latent space is smooth and continuous, meaning that nearby points in the latent space correspond to similar configurations in the original space. Smoothness of the latent space makes this easier to explore, interpolate over, and optimize. For example, in a process design, instead of working directly with a decision space that handles discrete variables (e.g., number of stages and the location of the feed stage), the VAE transforms these into a continuous latent representation. The smoothness of the latent space facilitates the development of surrogate models, such as GP models.  

\begin{figure}[!htp]
\centering
    \includegraphics[width=\textwidth]{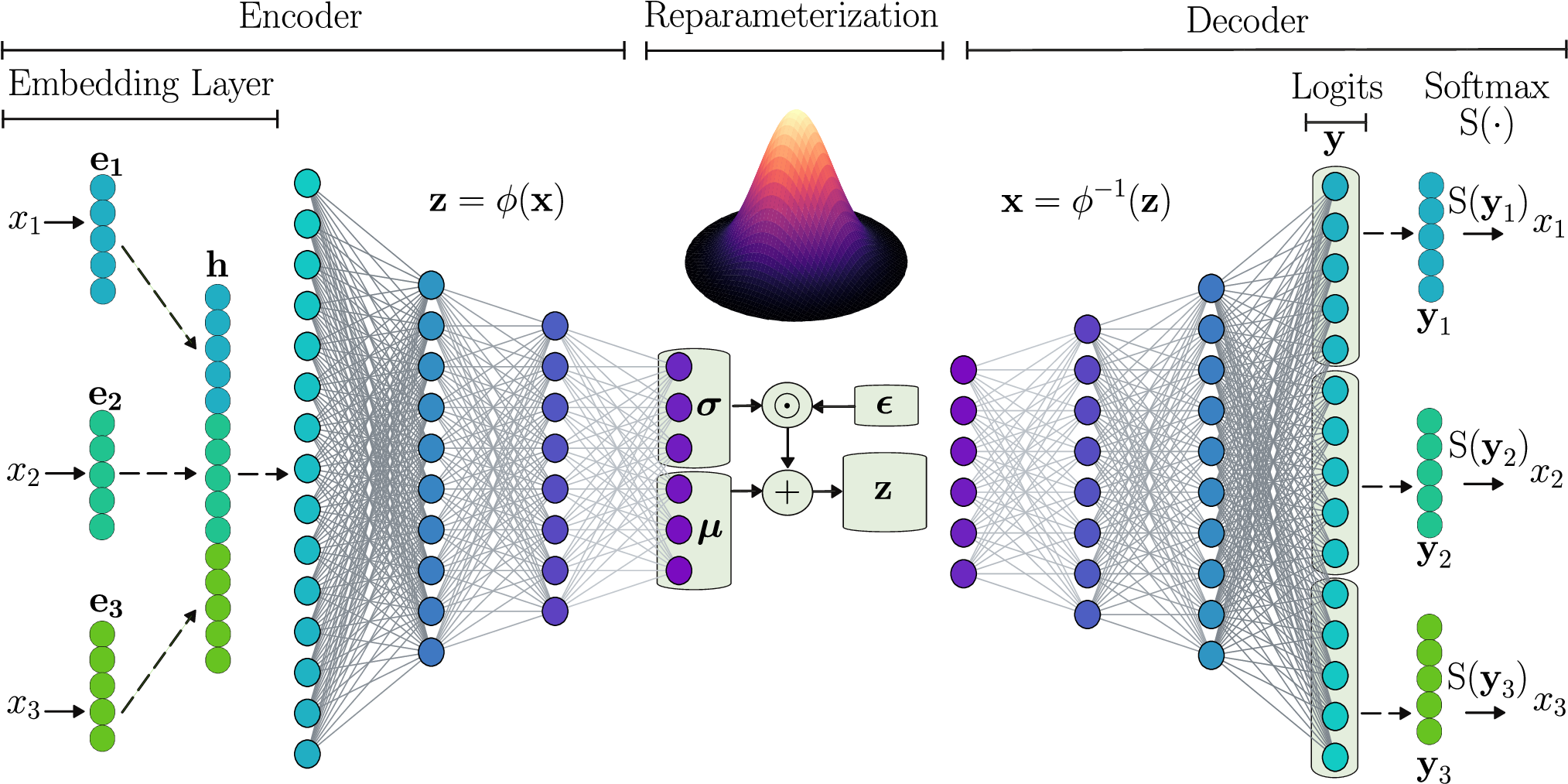}
    \caption{VAE architecture designed for a set of 3 discrete inputs $\mathbf{x} = [x_1, x_2, x_3]\in \mathcal{X}$; each discrete input is mapped into a continuous representation $\mathbf{e}$ through an embedding layer. The embeddings are concatenated into a vector $\textbf{h}$ and processed by the encoder. The encoder produces the mean vector $\boldsymbol{\mu}$ and the standard deviation vector $\boldsymbol{\sigma}$, which are used to generate the latent space $\mathcal{Z}$. The decoder maps points in the latent space $\mathbf{z}\in \mathcal{Z}$ back into the original space by producing logits $\mathbf{y}$ for each input variable; a softmax function is used to obtain probability distributions over the possible values of $x_1, x_2$, and $x_3$ to reconstruct the original discrete space (analogous to a classification model).}
    \label{VAE_structure}
\end{figure}

A schematic of a typical VAE architecture is presented in Figure \ref{VAE_structure}. Data for each input variable is transformed into a continuous vector $\textbf{e}$ using an embedding layer. Embedding enables learning dense, low-dimensional representations that are adjusted during the training process of the encoder network to capture semantic relationships between the embedded items, which enables the conversion of discrete input data into continuous, learnable representations, which are essential for many machine learning models. The embedded representations $\textbf{e}$ are combined into a single vector $\textbf{h}$ through concatenation. Subsequently, $\textbf{h}$ passes through a multilayer perceptron layer (MLP) to decrease the size until a given desired latent space dimension is reached. The embedding layer and the MLP are denoted as the encoder and are represented by the mapping $\phi$. The last layers of the encoder are used to perform a reparameterization to add stochasticity into the latent space $\mathbf{z}$ \citep{kingma2022}:
\begin{equation}
    \mathbf{z} = \boldsymbol{\mu} + \boldsymbol{\epsilon} \cdot \exp\left( \frac{1}{2} \log \boldsymbol{\sigma}^2 \right), \quad \boldsymbol{\epsilon } \sim \mathcal{N}(0,\mathbf{I})
\end{equation}
where $\boldsymbol{\epsilon}$ is standard Gaussian noise with mean zero and variance of one, $\boldsymbol{\mu}$ is the mean vector which are mean values for the latent variables and $\boldsymbol{\sigma}$ the standard deviation vector (that implicitly captures the range of possible values for the latent variables). Symbol $\mathbf{I}$ is the identity matrix, indicating that all components of $\boldsymbol{\epsilon}$ are independent standard normal variables. The probabilistic reparameterization enables the use of gradient-based optimization. The decoder (denoted as the inverse mapping $\phi^{−1}$), takes a point in the latent space $\mathbf{z}$ and predicts logits values $\mathbf{y}$, which are used to predict probabilities using the softmax function $\mathrm{S(\cdot)}$. This framework allows the reconstruction of the original input by categorical probabilities or by discrete class predictions.

The VAE is trained (the parameters of the encoder and decoder networks are learned) to minimize a composite loss function that captures the cross-entropy $\mathcal{L}_{\text{rec}}$ loss and the Kullback-Leibler divergence $\mathcal{L}_{\text{KL}}$ loss. The cross entropy $\mathcal{L}_{\text{rec}}$ captures individual entropies for each discrete variable and measures the ability to reconstruct the original space from the latent space (which is a measure of the reconstruction error): 
\begin{equation}
    \mathcal{L}_{\text{rec}} = - \sum_{n=1}^{N} \sum_{c=1}^{C} \text{log}\left(\frac{ e^{\mathbf{y}_{n,x_{n}}}}{\sum_{i=1}^{C} e^{\mathbf{y}_{n,i}}}  \right). 
\end{equation}
Here, $C$ is the number of possible classes/values for each discrete variable, and $N$ the total number of input discrete variables. Symbol $x_n$ denotes the index value used to measure the distance between the true class (discrete value) and the predicted class (reconstructed value). The Kullback-Leibler divergence $\mathcal{L}_{\text{KL}}$ measures how far the learned latent distribution $q({\mathbf{z}}|\mathbf{x}) \sim \mathcal{N}(\boldsymbol\mu, {\boldsymbol\sigma}^2)$ is from the target standard normal distribution $p(\mathbf{z}) \sim \mathcal{N}(0,\mathbf{I})$:
\begin{equation}
    \mathcal{L}_{\text{KL-D}} = \frac{1}{2} \sum_{j=1}^{D} \left( \exp(\log\boldsymbol\sigma_j^2) + \boldsymbol\mu_j^2 - \log\boldsymbol\sigma_j^2 -1\right)
\end{equation}
We note that the Kullback-Leibler divergence acts as a regularization term that aims to induce a structure on the latent space. 

The composite loss function is given by: 
\begin{equation}
    \mathcal{L} = \mathcal{L}_{\text{rec}} + \beta \cdot \mathcal{L}_{\text{KL}}
\end{equation}
where $\beta \in \mathbb{R}_{+}$ is a weighting parameter that trades-off the loss functions. 

In summary, the use of a VAE offers a couple of important benefits: i) it will represent discrete and continuous variables in a smooth and continuous space and ii) it will provide a bridge between the original problem formulation with advanced optimization tools (such as BO), which benefit from smooth search spaces. 

To illustrate the applicability of the VAE, we consider a simple problem of designing a reactor with the goal of maximizing the production of component $C$, which is produced simultaneously through a couple of competing routes, a consecutive pathway \((A \xrightarrow{k_1} B \xrightarrow{k_2} C)\) (referred to as Pathway 1) and a parallel direct pathway \((A \xrightarrow{k_3} C \xrightarrow{k_4} D)\) (referred to as Pathway 2). The kinetics of the reactor system are given by a set of elementary reactions, each characterized by a rate constant dependent of temperature $k_i(T)$ that follows the Arrhenius equation, \(k_i(T) = A_i \exp\left(-\frac{E_{i}}{RT}\right)\). We note from Pathway 2 that accelerating the formation of $C$ simultaneously enhances its undesired degradation to $D$. Each pathway has its own maximum production, which depends on both temperature and time. The goal is thus to determine the operating conditions that balance these competing effects and maximize the production of specie $C$, considering both pathways. The selection of the reaction pathway is a discrete design decision; this could be manipulated, for instance, by changing the catalyst.  The continuous decision variables include reaction time and temperature, which are discretized from their continuous spaces. 

\begin{figure}[!htp]
    \centering
    \includegraphics[width=1.0\textwidth]{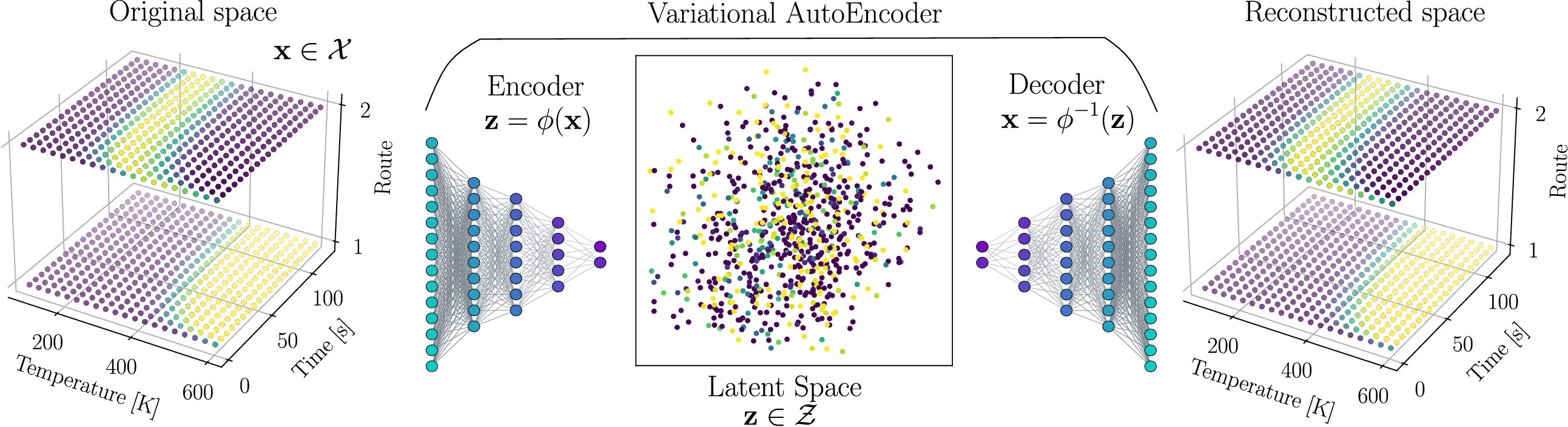}
    
    \caption{Schematic representation of the original discrete space $\mathcal{X}$ mapped through the decoder $\phi$ to obtain the latent continuous space $\mathcal{Z}$ from which we can reconstruct the original space using the decoder $\phi^{-1}$.}
     \label{toy_figure1}
\end{figure}

Within the discrete decision space, we evaluate the production of component $C$ along each pathway, and we can visualize the production of $C$ across the entire design space in Figure \ref{toy_figure1}, where each point on this grid corresponds to a combination of temperature, reaction time, and pathway. This visualization highlights the presence of a couple of local productivity maxima. The discrete space is projected to the continuous latent space $\mathcal{Z}$. By learning the encoder $\phi$ and decoder $\phi^{-1}$ mappings, a low-dimensional latent representation of the design space that captures its essential features is obtained; in addition, we obtain a  mapping to reconstruct the original variables from this representation. This approach will allow the exploration of feasible solutions in the latent space while ensuring valid operating choices for temperature, time, and pathways. 

For a suitable VAE architecture, one expect that the reconstruct of the original space with small reconstruction errors. A key challenge in enabling this is to find a suitable dimension of the latent space that minimizes the reconstruction loss. In Figure \ref{toy_figure2}, we illustrate the influence of the dimension of the latent space on the loss function. We note that the total (composite) loss stabilizes at a latent space dimension of eight. Here, we also show the profiles for each element of the loss as the VAE training progresses. Once the VAE can successfully reconstruct the original space from a latent space, we can proceed to use the latent space to search for optimal operating conditions.  In the next section we will discuss how this can be done systematically using BO.

\begin{figure}[!htp]
\centering
    \begin{subfigure}[b]{0.48\textwidth}
        \centering
        \includegraphics[width=1.0\textwidth]{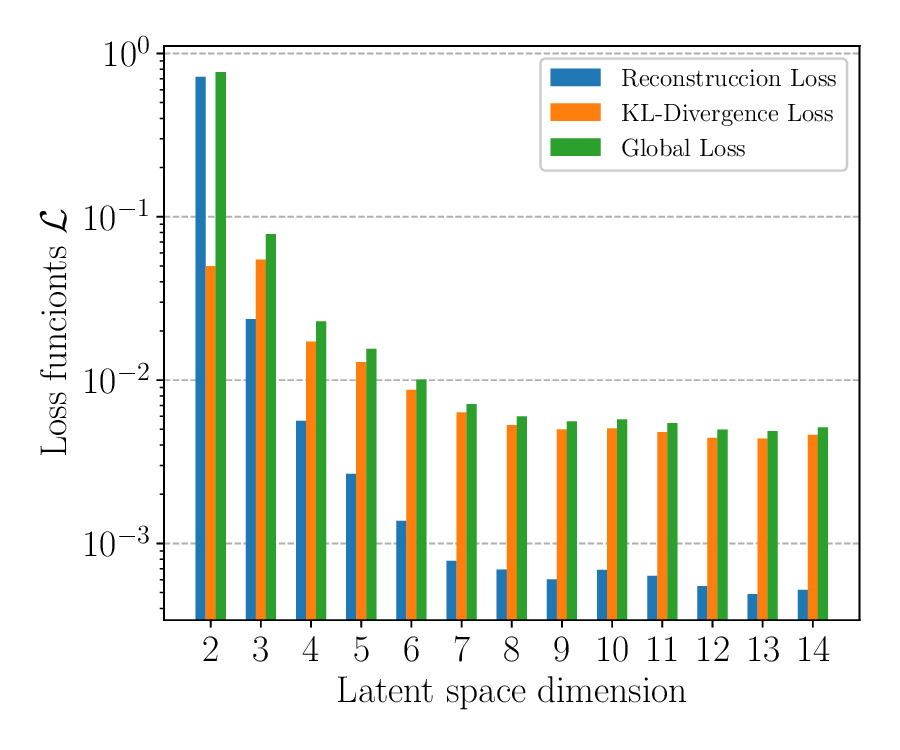}
        \label{toy_figure1_2}
     \end{subfigure}
    \begin{subfigure}[b]{0.48\textwidth}
        \centering
        \includegraphics[width=1.0\textwidth]{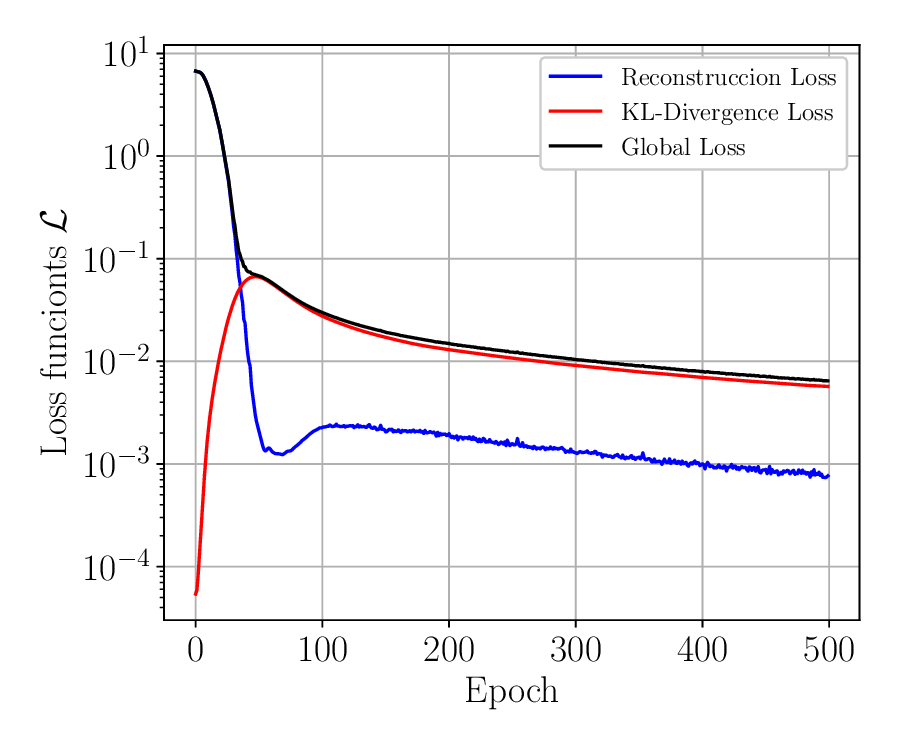}
        \label{toy_figure1_3}
     \end{subfigure}
     \vspace{-0.3in}
    \caption{Impact of the latent space dimension on the composite loss function $\mathcal{L}$ and its components: the reconstruction $\mathcal{L}_{\mathrm{rec}}$ and the KL-divergence $\mathcal{L}_{\mathrm{KL-D}}$ functions (left). Progression of loss functions during the training of a VAE with a latent dimension of eight.}
     \label{toy_figure2_hyperparameter}
\end{figure}

\section{Bayesian Optimization (BO) in VAE Latent Space}

Our goal is to solve the optimization problem:
\begin{equation}
    \min_{\mathbf{x} \in \mathcal{X}}\;\; f(\mathbf{x})
    \label{eq:bo_mixed}
\end{equation}
The key challenge is that the evaluation of the objective (cost) function $ f(\mathbf{x})$ at a point $\mathbf{x}\in\mathcal{X}$ requires a potentially expensive simulation. By transforming the original discrete design space $\mathcal{X}$ into a continuous space $\mathcal{Z}$ using a VAE, we can pose the optimization problem as:
\begin{equation}
    \min_{\mathbf{z} \in \mathcal{Z}} f(\phi^{-1}(\mathbf{z}))
    \label{eq:bo_mixed_modified}
\end{equation}
By changing the search space, we now aim to navigate the latent space so that the objective function is minimized; every point in the latent space is converted into a point in the discrete space. To navigate through the latent space, we will use a BO approach; this will build a surrogate model $m(\mathbf{z})$ that aims to approximate the composite mapping $f(\phi^{-1}(\mathbf{z}))$ at a given point in the latent space $\mathbf{z}\in \mathcal{Z}$. With the use of the surrogate model, the BO approach navigates through the original design space implicitly by navigating through the latent space.
\\

\begin{figure}[!htp]
\centering
    \includegraphics[width=1.0\textwidth]{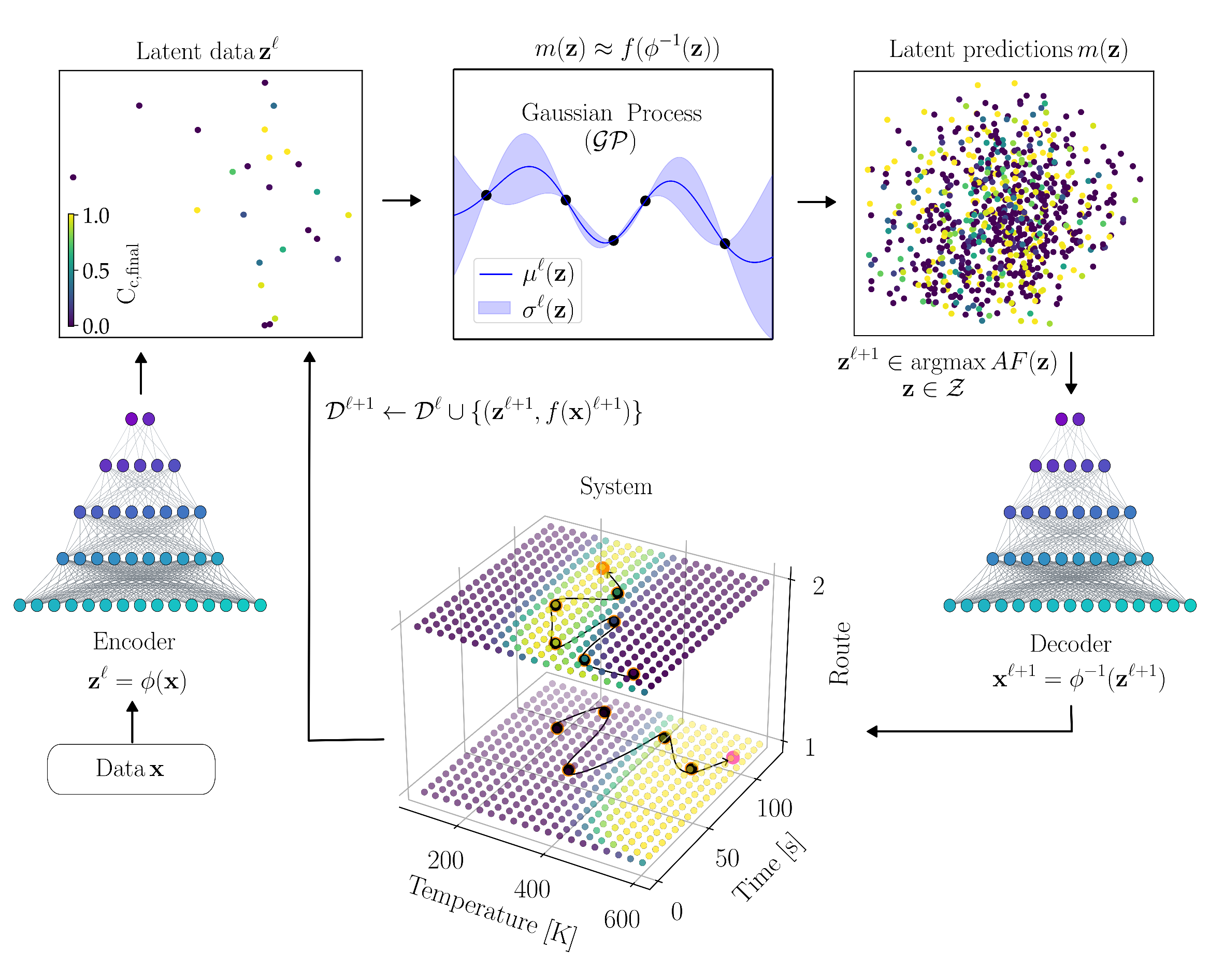}
    \vspace{-0.2in}
    \caption{Illustration of integrated BO+VAE framework. Data $\textbf{x}$ is mapped to $\textbf{z}$ using $\phi$. A surrogate model over the latent space $m(\mathbf{z})$ is built and iteratively selects promising points in the latent space $\mathbf{z}^{{\ell+1}}$ by optimizing an acquisition function ($AF$). Latent samples $\mathbf{z}^{\ell+1}$ are decoded to obtain a simulation at the new point $f(\phi^{-1}(\mathbf{z}^{\ell+1}))$. The surrogate model is refined by updating the dataset $\mathcal{D}^{\ell}$.}
    \label{toy_figure2}
\end{figure}

The combined VAE+BO framework is sketched in Figure \ref{toy_figure2}. At iteration $\ell$, we have a set of data points $\mathbf{x}^{\ell}\in\mathcal{X}$ with their corresponding costs $f(\mathbf{x}^{\ell})$. The input data points are encoded using $\phi$ to obtain the latent points $\mathbf{z}^\ell\in\mathcal{Z}$ and associated cost values $f(\phi^{-1}(\mathbf{z}^\ell))$. We build a probabilistic surrogate model $m^{\ell}(\mathbf{z})$ to approximate the composite mapping $f(\phi^{-1}(\mathbf{z}))$. The surrogate $m^\ell(\mathbf{z})$ is trained using the available cumulative database $\mathcal{D}_\ell = \{ \mathbf{z}^{k}, f( \phi^{-1}(\mathbf{z}^{k}))\}^\ell_{k=1}$. A Gaussian Process model ($\mathcal{GP}$) is used as a surrogate model, which is a probabilistic model that provides a mean prediction $\mu_f^\ell (\mathbf{z})$  of the cost function and the  uncertainty $\sigma_f^\ell (\mathbf{z})$ of the prediction.

We note that the VAE+BO framework includes a couple of approximation errors: the VAE reconstruction error (obtained by mapping from the latent to the original space) and the error of the surrogate model. An important observation is that, even if there are reconstruction errors, the refinement of the GP model should be able to correct for any discrepancies. However, we expect that the reconstruction error needs to be sufficiently small to ensure that points in the latent space are true reflections of points in the original space to have a consistent search. 

BO selects a new decision by optimizing an acquisition function (${AF}$) that balances exploration and exploitation. Exploration is driven mainly by the uncertainty, $\sigma_f^\ell (\mathbf{z})$ where new regions not explored are prioritized; exploitation focuses on the mean cost $\mu_f^\ell (\mathbf{z})$ predicted to refine the solution in a specific design space. To balance exploration and exploitation, we use the Expected Improvement ($EI$) as acquisition function. 

We are often interested in optimizing multiple objective functions to assess trade-offs. Here, we focus on bi-criteria problems of the form:

\begin{equation}
    \min_{\mathbf{z} \in \mathcal{Z}} (f_1(\phi^{-1}(\mathbf{z})),f_2(\phi^{-1}(\mathbf{z})) )
    \label{eq:MOBO_mixed}
\end{equation}

In this setting, the goal of BO is to discover the so-called Pareto frontier (set of non-dominated points); to do so, we use a simple procedure in which we build a GP model for each cost function and use the predictive means $\mu_f^\ell (\mathbf{z})$ to identify the predicted non-dominated points and thus predict the location of the Pareto front. At each point on the predicted Pareto frontier, we use the uncertainty information $\sigma_f^\ell (\mathbf{z})$ to determine the locations on the front that have the highest uncertainty. The next iteration $\mathbf{z}^{\ell+1}$ is selected by identifying the point on the Pareto front with the highest uncertainty; as we will illustrate below, this simple data acquisition approach can discover the Pareto frontier with few simulations. 

\section{Case Studies}

To assess the effectiveness of the proposed framework, we built a couple of detailed case studies that aim to design complex separation systems; the first system conducts liquid-liquid extraction of caprylic acid from water and the second system involves the separation of an azeotropic mixture by a dividing wall column. We implemented detailed simulation models for these systems using AspenPlus; our goal is to identify optimal designs for such systems by using as few expensive simulations as possible. 

\subsection{Liquid-liquid extraction of caprylic acid from water}

Caprylic acid (CA), a medium-chain fatty acid, is a platform chemical that has garnered significant attention due to its applications in the food and pharmaceutical industry \citep{ZHOU2019348, Shijian2015}. This component is mainly obtained through chain elongation of syngas during anaerobic fermentation. A challenge of this process is that fermentation yields highly diluted systems from which CA needs to be recovered; here, we design a liquid-liquid extraction system to address such challenge.
\\

\begin{figure}[!htp]
    \centering
    \includegraphics[width=0.7\textwidth]{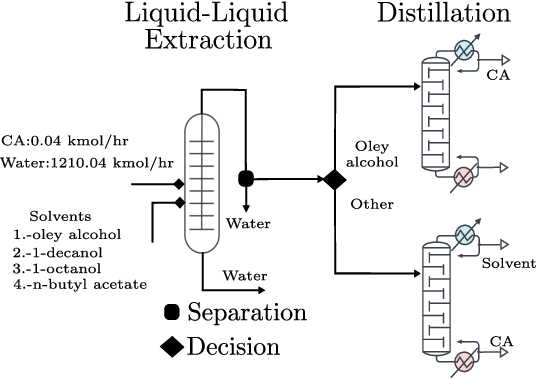}
    \caption{Process for the recovery of caprylic acid using liquid-liquid extraction. }
    \label{case1process}
\end{figure}

In \citet{Yuan2022}, an analysis was conducted to identify potential solvents (entrainers) that can be used to  capture CA with high selectivity. Among the eight solvents tested, they found that oleyl alcohol, 1-decanol,  1-octanol, and n-butyl acetate have high selectivity. A schematic of the process is shown in Figure \ref{case1process}. A bioreactor where the anaerobic fermentation takes place is followed by a centrifuge operation that splits biomass, leaving a mixture of CA and water at 298.15 K with a total molar flow rate of 1210.04 kmol/hr and 0.04 kmol/hr for CA and water, respectively. Liquid-liquid extraction is performed using a 10-stage column, and this task can be conducted with either of the four solvents mentioned above. The remaining water is separated, and the final step consists of a distillation column where the distillate rate or bottom rate are fixed to 0.042 kmol/hr (the value depends on the solvent selected). All units involved are considered assumed to operate at atmospheric pressure. 

The design aims to identify optimal solvents and configurations for the column that maximizes the molar fraction of CA ($x_{CA}$), which is a metric of {\em separation efficiency}. The second design objective is to minimize the reboiler duty of the column, which is a metric of {\em energy efficiency}. In other words, the goal is to identify designs that balance product purity and energy demand. Thermodynamic data used are reported in the Supplementary Information (SI). The design variables are given in Table \ref{variablesca};  all continuous variables were discretized. 

\begin{table}[!htp]
\centering
\caption{Design variables of caprylic acid recovery. \label{variablesca}} 
\begin{threeparttable}
    \begin{tabular}{l|ccc}

    \textbf{Design variables} & \textbf{Units}  & \textbf{Range} & \textbf{Domain} \\
    \hline
    Solvent Selection ($S$) & [-]  &  $\{0,1,2,3\}$\tnote{*}  & Categorical \\
    Reflux Ratio ($R$) & [-]  & $(0.1-3.0)$ & Continuous   \\ 
    Solvent Flowrate ($F_{S}$)  & kmol/hr  & $(1.5-10.0)$ & Continuous  \\
    Feed Stage ($NF_{stage})$  & [-] & $\{1-60\}$ & Discrete  \\  
    Total Number of Stages ($NT_{stages}$)  & [-] & $\{20-60\}$ & Discrete 
    \end{tabular}
\begin{tablenotes}
    \item[*] Oleyl alcohol, 1-octanol, 1-decanol, and n-butyl acetate are encoded using categorical values of 0, 1, 2, and 3, respectively. 
\end{tablenotes}
\end{threeparttable}
\end{table}

\begin{figure}[!htp]
    \centering
    \begin{subfigure}[b]{1.0\textwidth}
        \centering
        \includegraphics[width=1.0\textwidth]{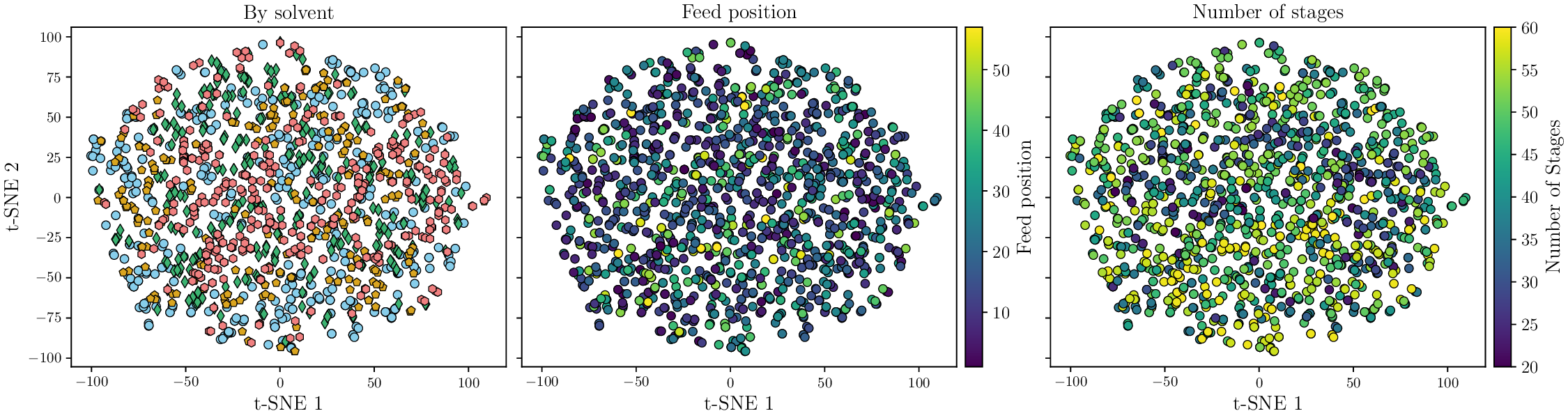}
        \caption{}
        \label{case1_alltsn}
    \end{subfigure}
    
    \begin{subfigure}[b]{0.45\textwidth}
        \centering
        \includegraphics[width=1.0\textwidth]{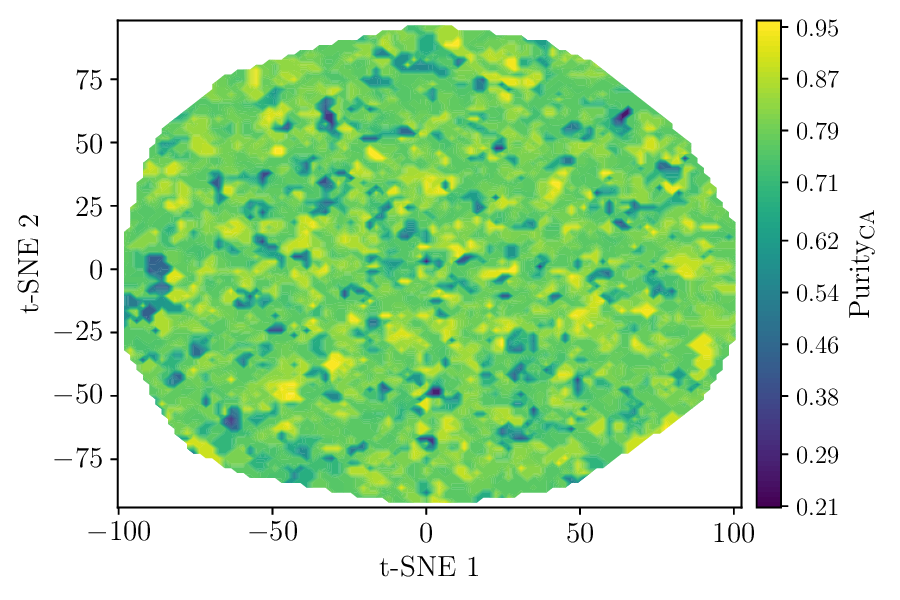}
        \caption{}
        \label{case1purity}
    \end{subfigure}
    \begin{subfigure}[b]{0.45\textwidth}
        \centering
        \includegraphics[width=1.0\textwidth]{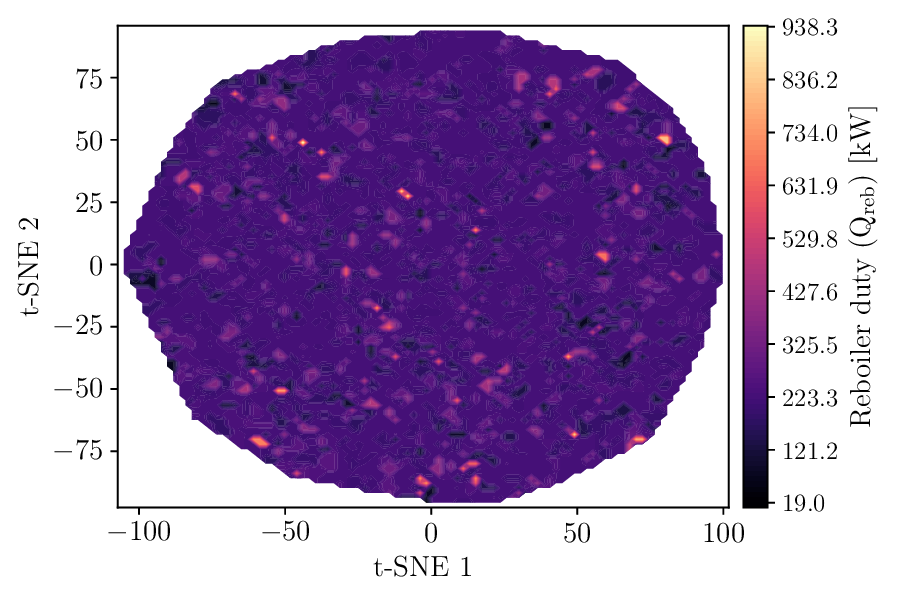}
        \caption{}
        \label{case1reboiler}
    \end{subfigure}
    \caption{Projection of design in the original space $\mathcal{X}$ into latent space $\mathcal{Z}$ using reduced with t-SNE; we highlight points corresponding to the solvent used, feeding position, and number of stages (a). Solvents are represented as follows: \legendnode{mskyblue}{circle} oley alcohol,\; \legendnode{mseagreen}{diamond} 1-decanol,\; \legendnode{mgoldenrod}{regular polygon, regular polygon sides=5} 1-octanol,\; \legendnode{mlightcoral}{regular polygon, regular polygon sides=6} n-butyl acetate. Visualization of the latent space in terms of values of the objectives of purity (b) and reboiler duty (c).}
    \label{results1}
\end{figure}

Figure \ref{case1_alltsn} visualizes the latent space $\mathcal{Z}$ in 2 dimensions; this visualization was obtained using t-distributed Stochastic Neighbor Embedding (t-SNE). Here, we highlight the points in the latent space associated with different solvents, feed positions, and number of stages. It is interesting to observe that the VAE induces some mixing of the discrete design variables. In Figures \ref{case1purity} and \ref{case1reboiler}, one can observe the latent space highlighting locations, different purity and reboiler duty (the objectives). One can see that the VAE generates a latent space that is ellipsoidal, highlighting the fact that the latent space is structured as a Gaussian random variable. Also, there are diverse locations in the latent space (corresponding to different designs) that achieve high purity and low reboiler duty.  These visualizations of the latent space leveraged the use of the GP model to predict the objectives in different regions of the search space. 

\begin{figure}[!htp]
    \centering
        \begin{subfigure}[b]{0.8\textwidth}
        \centering
        \includegraphics[width=0.9\textwidth]{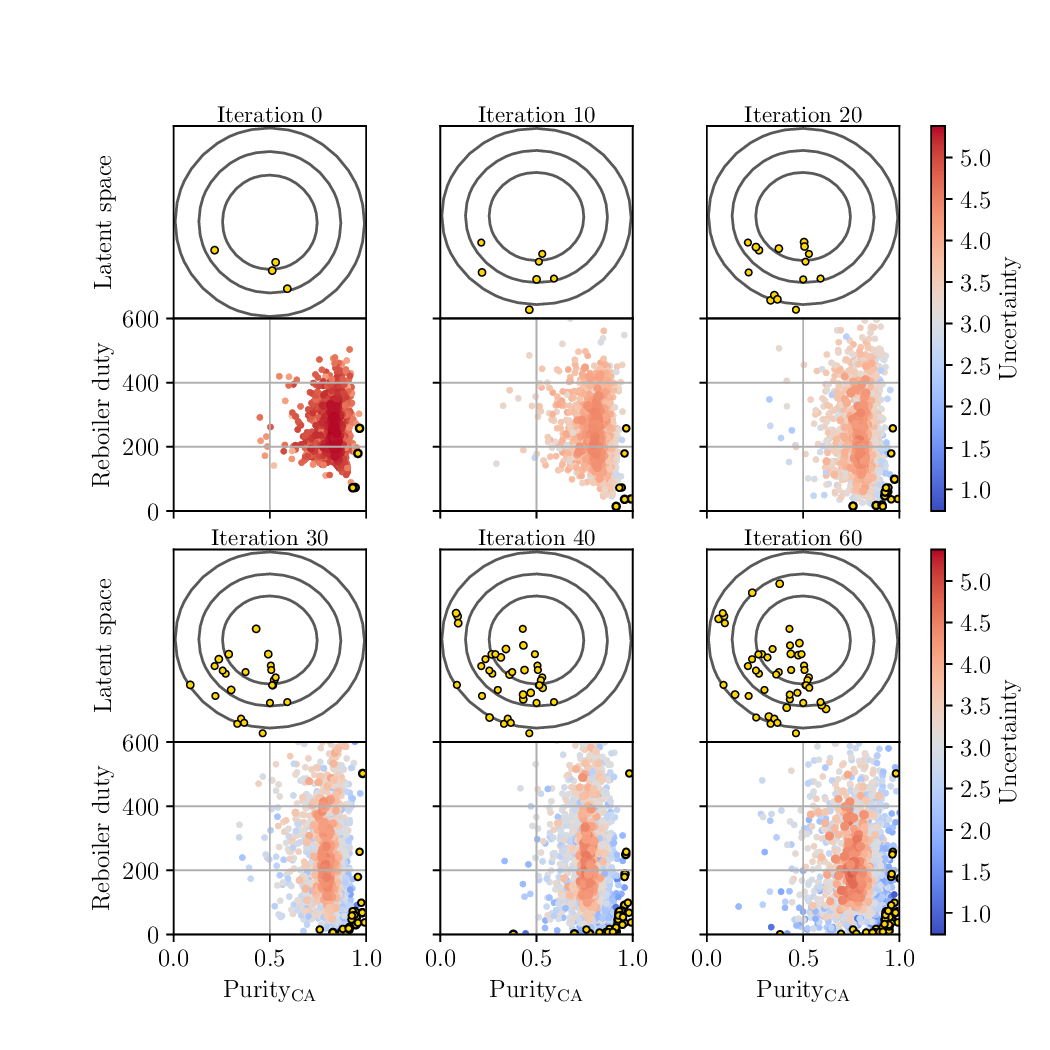}
        \caption{}
        \label{case1boprogress}
    \end{subfigure}
    
    \begin{subfigure}[b]{0.8\textwidth}
        \centering
        \includegraphics[width=1.0\textwidth]{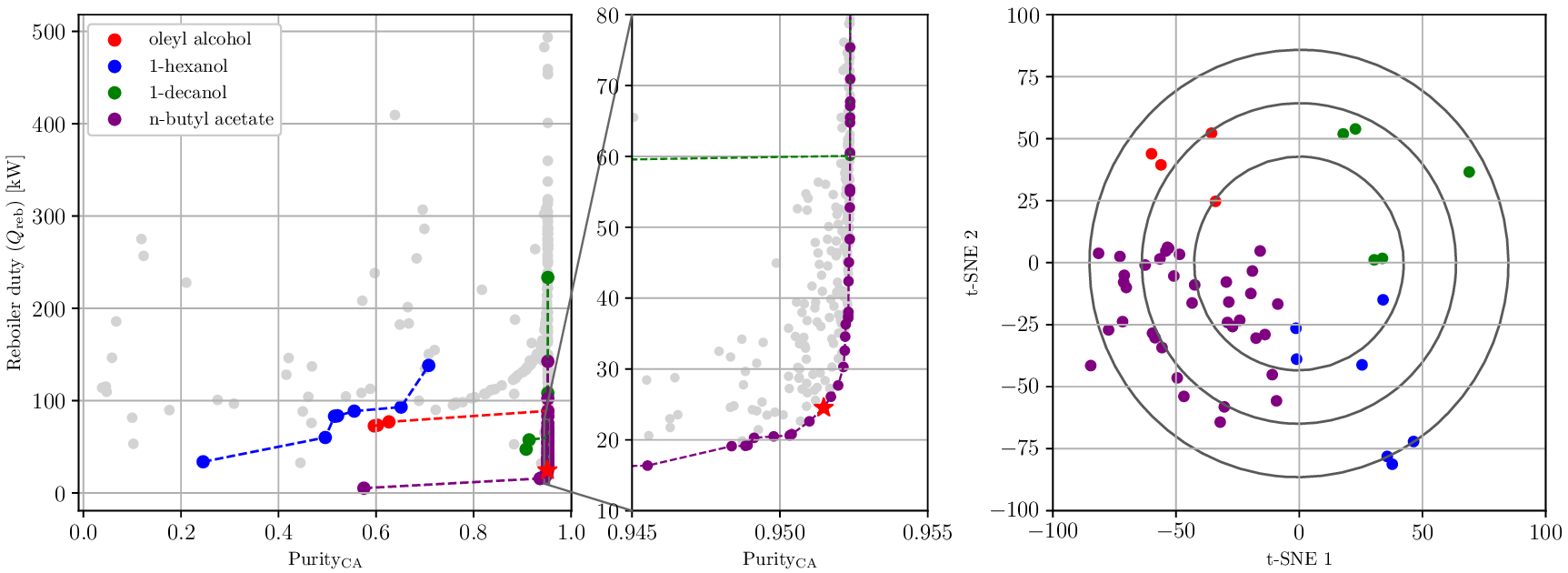}
        \caption{}
        \label{case1pareto}
    \end{subfigure}

    \caption{Evolution of the GP models and Pareto frontier along the BO search (a). Yellow dots represent the non-dominated (Pareto) points identified by the BO search. The concentric circles are level sets of the latent space, which highlight the positioning of the non-dominated points. We can see that the BO search quickly identifies a region in the space that corresponds to non-dominated points. Visualization of all simulations obtained with Aspen Plus and visualization of Pareto front highlighting different types of solvents and corresponding locations of the Pareto points in the latent space (b). We can see that the BO search identifies the Pareto frontier and we see that points along the Pareto front correspond to specific region of the latent space.}
    \label{results2}
\end{figure}

Figure \ref{case1boprogress} visualizes the evolution of the BO search over multiple iterations. The GP models were initially trained by using data collected from 50 initial simulations. In the first iteration, we can see that predictions have high uncertainty, which is attributed to the limited amount of data available. We can also see that most points are dominated points (not in the Pareto frontier). As more data are collected, the GPs are refined, allowing for the discovery of more non-dominated (Pareto) points and a reduction of uncertainty. For instance, in iteration 30, the uncertainty is significantly reduced along the Pareto frontier, and the predicted frontier shifts towards the bottom right. In iteration 60, we can see that most points begin concentrating on the Pareto frontier (indicated by low uncertainty), while the non-dominated region is not visited by BO (indicated by high uncertainty). This indicates that BO is strategically prioritizing the Pareto region.

All results from the simulations conducted with the AspenPlus simulator are shown in Figure \ref{case1pareto}. We confirm that, during the search, BO visits regions that are non-dominated points (not on the Pareto frontier). It is interesting to observe that different solvents tend to form different Pareto frontiers, but the best frontier is obtained for n-butyl acetate. By visualizing the Pareto points of this solvent in the latent space, one can observe that all these points cluster in the same region of the latent space. Here, we can also see that the points corresponding to other solvents are in other regions of the latent space. This highlights that the VAE is constructing a latent space with clearly defined regions. 

In Table \ref{tableresultados1}, a summary is given of different designs identified using BO, which enables the assessment of the diversity of the designs and associated objectives, thus providing a clear picture of the solution space explored in the search and the trade-offs involved. The average purity of compound CA ($x_{CA}$) was obtained as 93.6\%, with a relatively wide standard deviation of 8.3 \%, ranging from a minimum of 50.6\% up to a maximum of 95.2\%. This indicates that high purities can be achieved, although the system exhibits significant deviations; these deviations are correlated with the necessary amount of solvent required to carry out the separation. The average solvent feed flow rate $F_{S}$ was 3.11 kmol/hr, with a standard deviation of 1.03 kmol/hr, which ranged from 1.80 kmol/hr to 6.15 kmol/hr, reflecting the trade-off between solvent addition and separation performance. We can see that the compromise solution selects n-butyl acetate, achieves a purity of 95.14 \%, and has a reboiler duty of 24 kW. 

\begin{table}[!htp]
\centering
\caption{Descriptive statistics of non-dominated points, dominated for the caprylic acid recovery problem and configuration for the compromise solution found.}
\label{tab:stats_ca}
\begin{threeparttable}
\begin{tabular}{l|cccccccc}
\hline
    Point &  Metric & $x_{CA}$ & $Q_{reb}$ & $R$ & $F_{S}$ & $S$\tnote{*} & $NF_{stage}$ & $NT_{stages}$ \\
        &          & [\%]   & [kW]      & [-] & [kmol/hr] & [-] & [-] & [-] \\
\hline
\multirow{4}{*}{Non-Dominated} & mean & 93.60 & 51.0  & 0.400 & 3.108 & 3 & 24 & 45 \\
                               &  std   & 8.3 & 33.7  & 0.277 & 1.032 & 1 & 13 & 11 \\
                               &  min   & 50.6 & 10.47 & 0.110 & 1.799 & 2 & 1  & 20 \\
                               & max   & 95.2 & 148.29 & 1.168 & 6.145 & 3 & 57 & 60 \\
\hline
\multirow{4}{*}{Dominated}  & mean & 88.92 & 112.9 & 1.260 & 3.235 & 2 & 23 & 44 \\
                            &  std & 17.44 & 90.4  & 0.712 & 1.506 & 1 & 14 & 11 \\
                            &  min & 3.60  & 16.6  & 1.112 & 1.112 & 0 & 1  & 20 \\
                            & max  & 95.2  & 493.6 & 2.999 & 9.908 & 3 & 59 & 60 \\
\hline
 Compromise &   &  95.14 & 24.54 & 0.1551 & 2.25 & 3 & 29 & 44 \\  
\hline
\end{tabular}
\begin{tablenotes}
        \item[*]  Solvents oleyl alcohol, 1-octanol, 1-decanol, and n-butyl acetate are encoded with the categorical values 0, 1, 2, and 3.
\end{tablenotes}
\end{threeparttable}
\label{tableresultados1}
\end{table}

For the solvent type ($S$), 1-decanol and n-butyl acetate solvents were found to be the best options in terms of energy requirements. In terms of the reflux ratio ($R$), the average value was 0.40, with a minimum of 0.11 and a maximum of 1.17, indicating that both highly economical and more energy-intensive operating modes were explored. The mean reboiler duty (\(Q_{reb}\)) found was 51.0 kW, ranging from 148.2 to 10.5 kW with a standard deviation of 33.7 kW in the most energy-efficient cases. The feed stage number ($NF_{stage}$) had an average of 24 stages with a standard deviation of 13 stages, extending from a minimum of 1 to a maximum of 57. In contrast, the total number of stages ($NT_{stages}$) averaged 45, with a minimum of 20 and a maximum of 60. This diversity emphasizes the significant impact of operating conditions on the purity and energy requirements of the distillation process. For the dominated points, the average purity decreases to 88.92, but also provides designs with the maximum recovery. Regarding the reflux ratio, the values of the dominated points have a higher average and reach the upper exploration limit of 2.999. The same behavior can be observed in the case of solvent flow, where the average corresponds to 3.235 kmol/hr, with the maximum value explored being 9.908 kmol/hr. These dominated points exhibit a wide variety of designs for the different solvents tested with the average provided by 1-decanol as the predominant solvent, which indicates that the BO search is exploring a broad range of designs along the search.

Our findings emphasize the impact of solvent selection on the inherent trade-offs between energy efficiency and column design. The selection of the best solvent can also be justified by comparing the vapor-liquid equilibrium diagrams for each entrainer (see the SI for equilibrium diagrams). From all solvents tested, n-butyl acetate has the lowest boiling point of 463 K; oley alcohol, 1-decanol, and 1-octanol have boiling points of 623 K, 503 K, and 468 K, respectively. Nonetheless, the boiling points of n-butyl acetate and 1-octanol are close, making their selection non-trivial and thus requiring an optimization procedure.

\subsection{Dividing Wall Column}

The incorporation of thermally coupled units in chemical processes has generated significant attention due to their benefits by reducing capital and auxiliary utilities costs \citep{Czarnecki3c00302,TANG2025130030,HERNANDEZ19991005}. Dividing wall columns (DWC), in particular, have been applied in various processes worldwide \citep{GUTIERREZGUERRA2009145}. These units have been used for a wide range of applications, including difficult separations (e.g., mixtures with azeotropes). An application we consider here is the separation of organic mixtures generated in the pharmaceutical industry, specifically mixtures of methanol and dichloromethane (DCM). The design of DWCs is challenging due to the tight integration of the system components.

\begin{figure}[!htp]
    \centering
    \includegraphics[width=\textwidth]{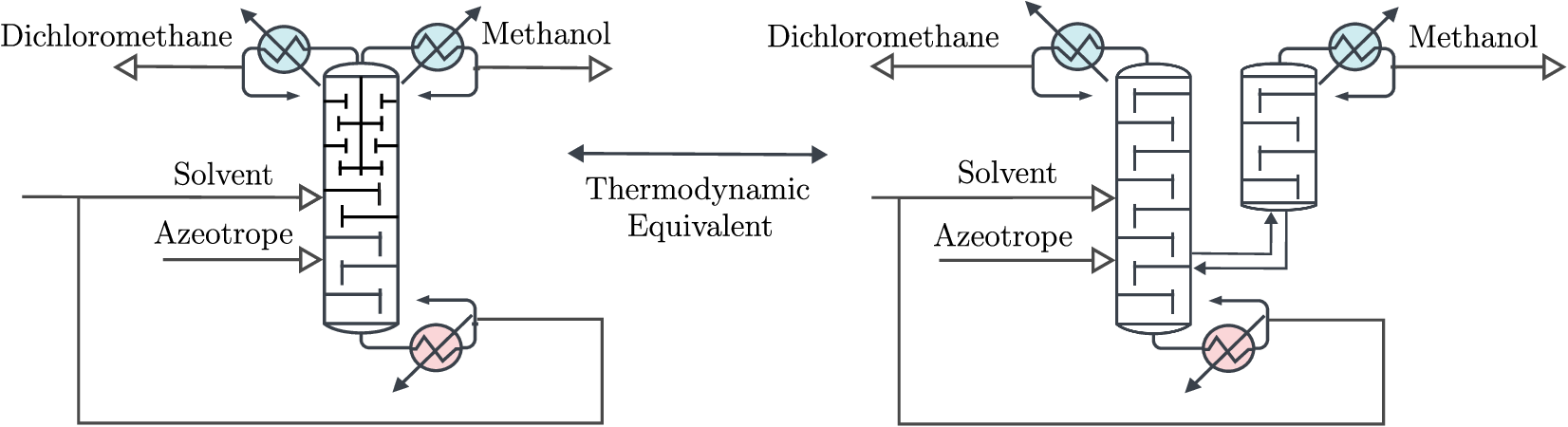}
    \caption{Thermodynamic representation of a dividing wall column for the separation of an azeotropic mixture of dichloromethane and methanol.}
    \label{case2process}
\end{figure}

A reported work has focused on identifying designs of DWCs that lead to high purity DCM ($\geq$ 99.90\%) \citep{Jing2021}. We extend here such work by evaluating the inherent safety of these columns as part of the optimization problem. The metric used is the Fire and Explosion Damage Index (FEDI), which was proposed by \citet{KhanAbbasi} and has been used for the early risk evaluations \citep{Landero8b04146, ORTIZESPINOZA2021225}. The index is included in this work as an additional objective function of this optimization problem. It is assumed that a lower DCM composition indicates lower purity of methanol. The design variables available to optimize the DWC are listed in Table \ref{variablesdwc}. A feed flow rate of 100 kmol/h was considered with a couple of possible variations corresponding to equimolar and near-azeotropic concentrations. In the equimolar scenario, both distillate rates are set at 50.05 kmol/hr, while the azeotropic feed values are set to 83.07 and 16.97 kmol/hr. The thermodynamic equivalent intensified column is shown in Figure \ref{case2process}. For the AspenPlus simulation, the NRTL thermodynamic model was used. The information needed to compute FEDI values and classification ranges is provided in the supplementary information.

\begin{table}[!htp]
\centering
\caption{Design variables to optimize the DWC and their search ranges. \label{variablesdwc}} 
\begin{threeparttable}
    \begin{tabular}{l|ccc}

    \textbf{Design variables} & \textbf{Units}  & \textbf{Range} & \textbf{Domain} \\
    \hline
    Reflux ratio ($R_{DWC}$) & [-]  & $(2.0-5.0)$ & Continuous   \\ 
    Vapor flowrate ($V_{DWC}$)  & [kmol/hr] & $(50-100)$ & Continuous \\
    Solvent to Feed Ratio ($S/F$)  & [-]  & $(1.5-10.0)$ & Continuous  \\
    Solvent feed stage ($NF_{solv})$  & [-] & $\{1-50\}$ & Discrete  \\  
    Azeotrope feed stage ($NF_{azeo}$) & [-] & $\{1-50\}$ & Discrete  \\  
    DWC Trays ($DWC_{trays}$) & [-] & $\{1-50\}$ & Discrete \\
    Total stages ($NT_{DWC}$) & [-]  &  $\{30,50\}$  & Discrete \\
    \end{tabular}
    
\end{threeparttable}
\end{table}

\begin{figure}[!htp]
    \centering
    \includegraphics[width=\textwidth]{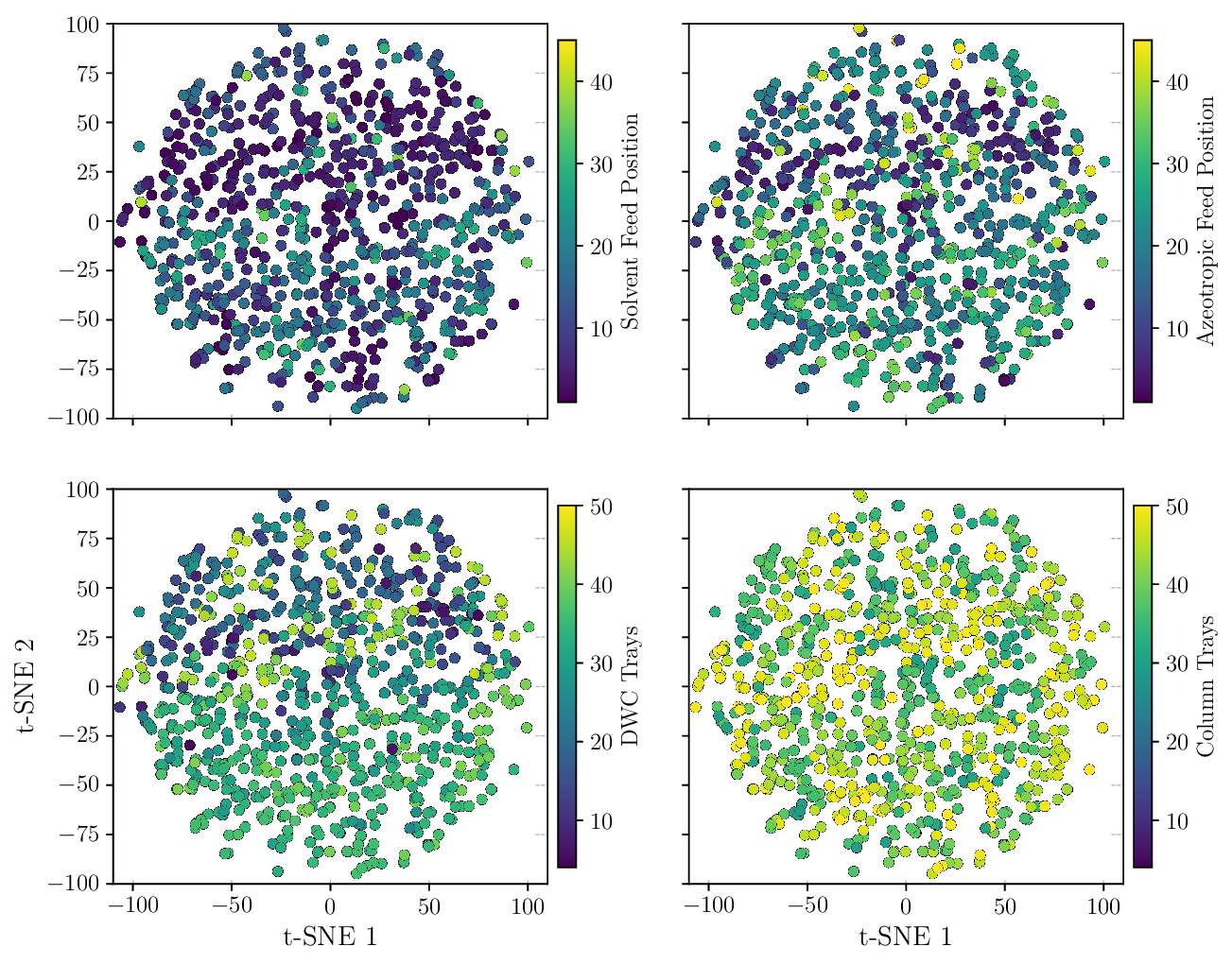}
    \caption{Projection of random sample points of the latent space $\mathbf{z}$ reduced with t-SNE and their corresponding distribution for the feeding positions of the solvent, azeotrope, size of diving wall, and total column trays.}
    \label{case2tsneall}
\end{figure}

As in the previous case study, we sampled points in the latent space and visualized them using t-SNE. Figure \ref{case2tsneall} shows this distribution with respect to the discrete variables: solvent and azeotrope feeding position, number of dividing wall stages, and total number of dividing column stages. We can see again that the VAE generates a mixed space that is continuous and follows ellipsoidal regions.

\begin{figure}[!htp]
    \centering
    \begin{subfigure}[b]{1.0\textwidth}
    \centering
        \includegraphics[width=0.8\textwidth]{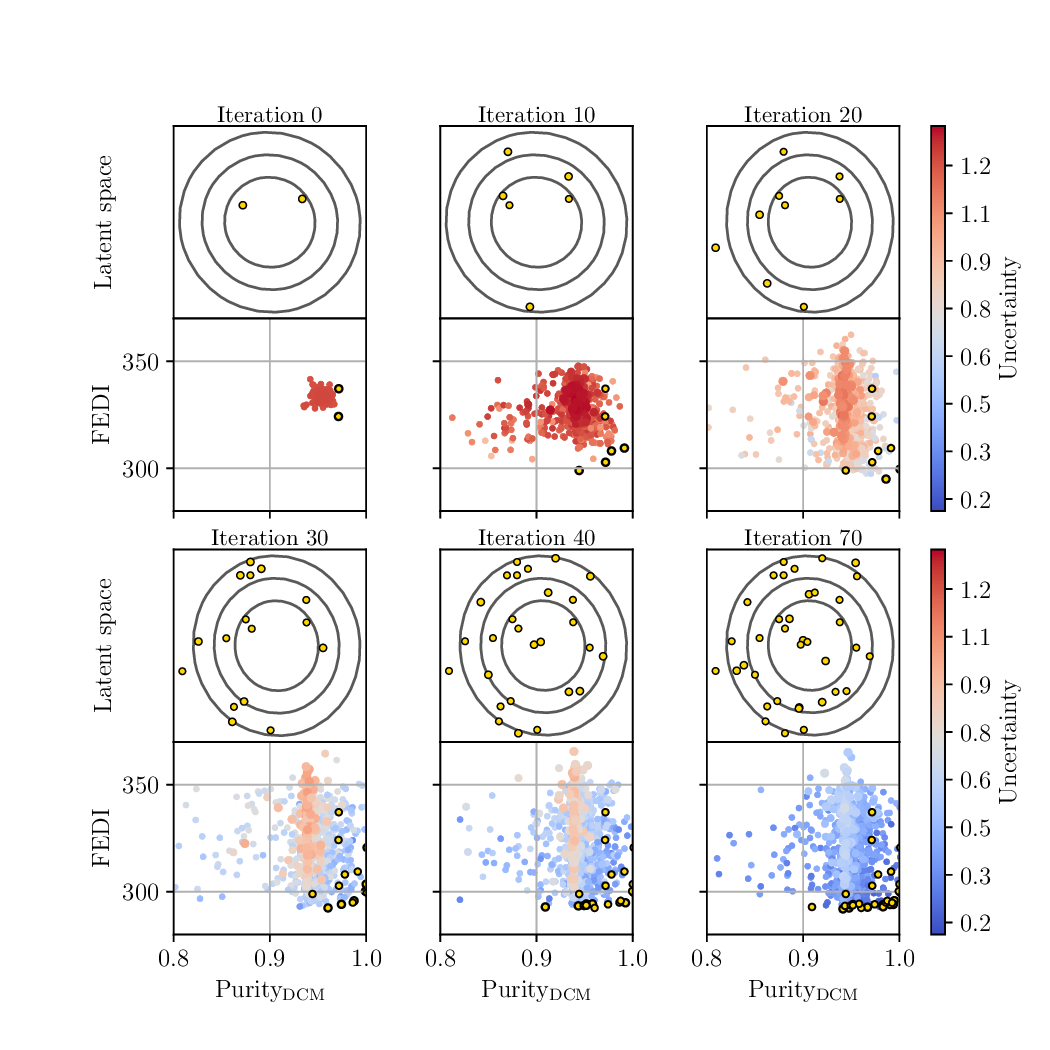}
        \caption{}
        \label{results3_1}
    \end{subfigure}
    
    \begin{subfigure}[b]{0.45\textwidth}
        \includegraphics[width=0.9\textwidth]{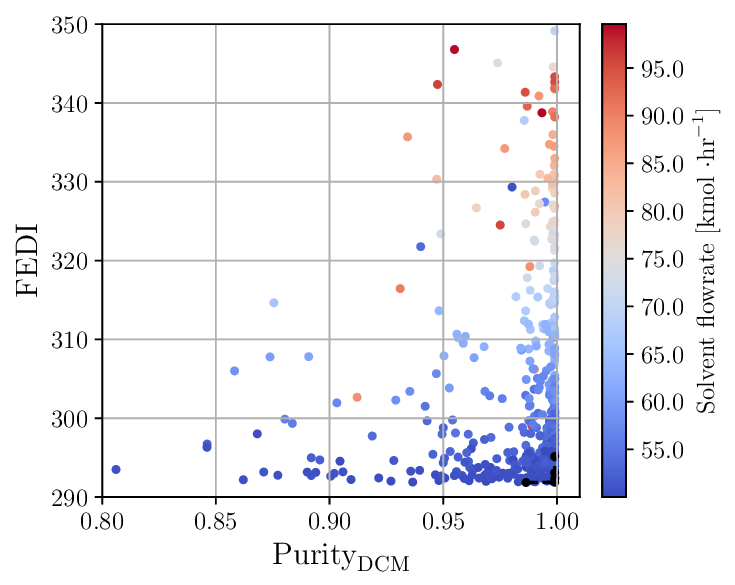}
        \caption{}
        \label{results3_2}
    \end{subfigure}
    \begin{subfigure}[b]{0.45\textwidth}
        \includegraphics[width=0.9\textwidth]{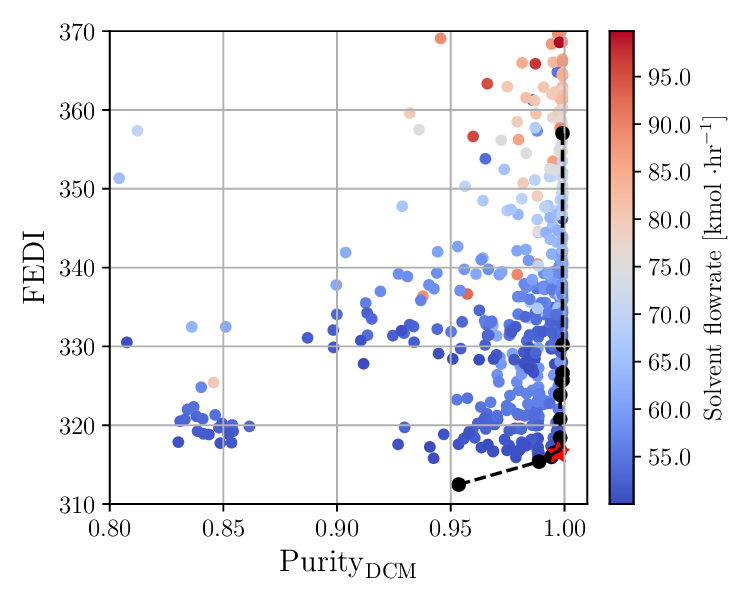}
        \caption{}
        \label{results3_3}
    \end{subfigure}
    \caption{ Evolution of BO search thought different iterations with their corresponding GP predictions (a). Yellow dots represent the non-dominated (Pareto) points collected. Feasible designs were found for the dividing wall column for equimolar composition (b) and azeotropic composition (c). Points are color-coded according to the solvent flow rate, which is a critical factor for the inherent safety metric. Black dots and dashed lines show the Pareto front.}
    \label{results3}
\end{figure}

The evolution of the BO search and their corresponding GP model predictions are shown in Figure \ref{results3_1}. Initially, with 50 data points, the predictions are highly restricted to the range of 0.9-1.0 mol fraction and 370-348 FEDI values. Once new feasible designs are discovered and collected, this space increases mainly in the FEDI region, where high uncertainty is located behind the Pareto front. From iteration 30 onward, it can be seen that the search intensifies within the area of interest. This progressive refinement is based on the scenario where an equimolar composition is assumed. Similarly, it can be seen that the latent points collected throughout the optimization process are distributed well throughout the space; this highlights that there is no particular region in the design space under which the Pareto solutions cluster (which means that this is a complicated search). 

Figure \ref{results3} also illustrates the relationship between composition and the inherent safety FEDI metric. Although the composition of the distillates is relaxed, the intensified unit is still classified as hazardous. In Figure \ref{results3_2}, the index FEDI varies within the range 294-350. However, when the composition of the azeotropic feed changes as shown in Figure \ref{results3_3}, these limits shift to 310-370. These results suggest that the inherent safety index is sensitive to composition perturbations or variations, and that multiple configurations can provide the purity specified, but in some cases, deteriorate the process safety. 

The designs identified in Figure \ref{results3} were classified into dominated and non-dominated solutions, and their corresponding statistics are shown in Table \ref{tableresultados2}. While the results of \citet{Jing2021} correspond to the solution of a mixture with a composition close to the azeotrope, these findings show that optimal solutions considering economy or safety considerations can have conflicts between them. According to their solution, the design that minimized costs operated with a reflux ratio of 0.25 in the main column and incorporated a dividing wall spanning 22 trays, resulting in a total of 31 trays for the DWC. Additionally, the necessary $S/F$ ratio to meet the purity requirements was 0.48, and the optimal steam split ratio was 0.81. Furthermore, when the objective function changed to the inherent safety, minimizing FEDI reduced this ratio substantially to 0.21 and increased the reflux ratio. Overall, our analysis indicates that the non-dominated solutions maintain consistent purity while exhibiting high variations in purity and safety performance, primarily due to adjustments in the reflux ratio and solvent-to-feed ratio. The dominated points show greater variability with respect to purities, reaching a minimum value of 57.43 \% and a maximum FEDI of 391. 

Statistics associated with column dimensionality show that both sets of points explore similar regions, implying that these variables do not have a significant impact on the safety indicator. Therefore, the critical factor identified for these distillation configurations is the amount of solvent used to achieve the separation. From a safety standpoint, minimizing the quantity of flammable materials processed within the unit significantly enhances the overall safety of the unit operation. With the proposed methodology, we effectively navigated and evaluated a complex structure, which can be adapted to enhance more complicated systems, such as reactive dividing wall columns or electrified distillation columns.

\begin{table}[!htp]
\centering
\caption{Non-dominated and dominated designs for divided wall column.}
\begin{threeparttable}
\begin{tabular}{lccccccccc}
\hline
 & $x_D$  & FEDI & $R_{DWC}$  & $S/F$ & $V_{DWC}$ &  $NF_{solv}$ & $NF_{azeo}$ & $DWC_{trays}$ & $NT_{DWC}$ \\
 & [\%] & [-] & [-] & [-]  & [kmol/hr] & [-] & [-] & [-] & [-]  \\
\hline
\multicolumn{10}{c}{Non-Dominated Points} \\ %
\hline
mean & 99.45 & 314 & 2.96  & 0.56 & 54.95  & 7 & 24 & 33  & 43 \\
std  & 1.08  & 17  & 0.81  & 0.13 & 10.28  & 6 & 4 &  6 & 6 \\
min  & 95.35 & 291 & 2.13  & 0.50 & 50.11  & 2 & 20 & 26 & 35 \\
max  & 99.91 & 357 & 4.45  & 0.97 & 95.04  & 22 & 36 & 45 & 50 \\
\hline
\multicolumn{10}{c}{Dominated Points} \\ %
\hline
mean & 97.25 & 321 & 3.47 & 0.60 & 54.30  & 10 & 21 & 31  & 43 \\
std  & 5.03  & 23  & 0.71 & 0.12 & 7.32  & 7 & 8 & 8  & 5 \\
min  & 57.43 & 292 & 2.00 & 0.50 & 16.05  & 1 & 2 &  5 & 30 \\
max  & 99.91 & 391 & 4.99 & 0.99 & 99.60  & 40 & 43 &  46 & 50 \\
\hline
c.s\tnote{1} & 99.91 & 293.5 & 4.2 & 0.50 & 52.5 & 4 & 34 & 38 & 42 \\
c.s\tnote{2} & 99.72 & 316.1 & 2.8 & 0.50 & 53.1 & 3 & 18 & 31 & 50 \\
\hline
\end{tabular}
\begin{tablenotes}
\item[1] Optimal solution for an \textit{equimolar} feed.
\item[2] Optimal solution for an \textit{azeotropic} feed.
\end{tablenotes}

\end{threeparttable}
\label{tableresultados2}
\end{table}

\section{Conclusions and Future Work}


This work illustrates how to transform discrete design spaces into continuous one using variational autoencoders (VAEs). We show that this transformation facilitates the implementation of simulation-based optimization approaches, such as Bayesian Optimization (BO). The proposed framework has been implemented on a CSTR reactor where all decision variables were discretized. This shows the potential applicability to real-world scenarios where experimental campaigns are limited by cost or time, resulting in exploration limited by a fixed budget. In addition, we showed the applicability of the approach to the design of conventional and intensified distillation processes, where the numerical convergence and the number of evaluations represent a computational limitation. In this sense, this approach can be applied to more complex distillation columns, such as multi-product dividing wall columns and thermally coupled reactive columns, where exploration costs are even higher. We show that the proposed approach is effective at navigating the design space and can quickly discover optimal solutions. The performance of the proposed approach has been evaluated via computational experiments.

As part of future work, we would like to conduct additional studies to determine the ability of transforming different types of design spaces and to provide a more rigorous justification of performance. We are also interested in exploring the ability of the VAE approach to represent discrete spaces that have complex geometries (as identified by constraints) and we are interested in using the design space transformation approach to other types of problems, such as purely continuous optimization problems.  On the other hand, a critical challenge of this methodology lies in the structure of the VAE, such as the architecture of the encoder and decoder, and the hyperparameters used. Improvement of the hyperparameter is essential to reduce reconstruction loss, avoiding unfeasible reconstructed points. In addition, the dimensionality of the latent space plays a crucial role in determining the overall performance of the VAE. While it provides a compact representation of the underlying data, an excessively high or low latent dimension can limit the model’s expressiveness and generalization ability. To mitigate this constraint, the search space can be adaptively partitioned through Bayesian Optimization techniques combined with trust-region methods, thereby accelerating the discovery of optimal solutions and enhancing the construction of a well-defined Pareto front.

\section*{Acknowledgments}

The development of this project was supported by a SECIHTI (Mexico) scholarship granted to Gabriel Hernández-Morales and the National Science Foundation. 

\bibliographystyle{plainnat} 
\bibliography{references}


\end{document}